# Origin of the numerals
# Al-Biruni's testimony

Ahmed Boucenna
Laboratoire DAC, Department of Physics, Faculty of Sciences,
Ferhat Abbas University 19000 Stif, Algeria
aboucenna@wissal.dz    aboucenna@yahoo.com

**Abstract**
 The origin of the numerals that we inherited from the arabo-Islamic civilization remained one enigma. The hypothesis of the Indian origin remained, with controversies, without serious rival. It was the dominant hypothesis since more of one century. Its partisans found to it and constructed a lot of arguments. The testimonies of the medieval authors have been interpreted to its advantage. The opposite opinions have been dismissed and ignored. An amalgam between the history of our modern numerals and the Indian mathematics history is made. Rational contradictions often passed under silence.
A meticulous observation of the numerals permits to affirm that our numerals are in fact more or less modified Arabic letters. The "Ghubari" shape of the numerals shows that the symbol of a numeral corresponds to the Arabic letter whose numerical value is equal to this numeral. The numerals don't have a simple resemblance with some Arabic letters, but every number looks like the Arabic letter whose numerical value is equal to this numeral. The elements of the "Abjadi" calculation gives us a theoretical support, independent of the letters and numerals, witch explains our observation. Besides a re-lecture of the testimonies of the medieval authors, particularly the testimony of Al-Biruni, that is probably at the origin of all others testimonies speaking of the Indian origin of the numerals, is in agreement with the fact that our numerals are Arabic letters.
We have there a second way concerning the origin of our modern numerals that is only to its beginnings. The deepened researches are necessary to understand the history of our numerals better. A rigorous re-lecture of the medieval testimonies with a new mind imposes itself.

MCS : 01A30

***Keywords*** : "Ghubari" numeral, modern numeral, "Mashriki" numeral, Arabic letters, Hebrew letters, letter numerical value, "Abjadi" order,

**Résumé**
L'origine des chiffres que nous avons hérité de la civilisation arabo-islamique est restée une énigme. L'hypothèse de l'origine indienne est restée, avec des controverses, sans rival sérieux. Elle fut l'hypothèse dominante depuis plus d'un siècle. Ses partisans lui ont trouvé et construit beaucoup d'arguments. Les témoignages des auteurs médiévaux ont été interprétés à son avantage. Les avis contraire ont été purement et simplement écartés voire ignorés. Un amalgame entre l'histoire de nos chiffres modernes et l'histoire des mathématiques indiennes est fait. Les contradictions rationnelles sont souvent passées sous silence.
Une observation méticuleuse des chiffres permet d'affirmer que nos chiffres ne sont en fait que des lettres arabes plus ou moins modifiées. La forme "Ghubari" des chiffres montre que le symbole du chiffre correspond à la lettre arabe dont la valeur numérique est égale à ce chiffre. Les chiffres n'ont pas une simple ressemblance avec certaines lettres arabes, mais chaque chiffre ressemble à la lettre arabe dont la valeur numérique est égale à ce chiffre. Les éléments du calcul ''Abjadi'' nous donne un support théorique indépendant des lettres et des chiffres qui explique notre observation. De plus une relecture des témoignages des auteurs médiévaux, particulièrement le témoignage d'Al-Biruni, qui est probablement à l'origine de tous les autres témoignages parlant de l'origine indienne des chiffres, est en parfait accord avec le fait que nos chiffres sont des lettres arabes.

Nous avons là une deuxième voie concernant l'origine de nos chiffres modernes qui n'est qu'à ses débuts. Des recherches approfondies sont nécessaires pour mieux comprendre l'histoire de nos chiffres. Une relecture rigoureuse des témoignages médiévaux avec un nouvel esprit s'impose.

**1. Introduction**

Before approaching the question of the origin of the modern numerals, let us first begin by recalling the following terminology. In this work, we designate by :
- Modern numerals or Arabic modern numerals, the numerals whose symbols are the following: 0, 1, 2, 3, 4, 5, 6, 7, 8, 9, which are used in Sciences and Technologies.
- "Mashriki" numerals, or Arabic "Mashriki" numerals, whose symbols are:

٠ ١ ٢ ٣ ٤ ٥ ٦ ٧ ٨ ٩

which are currently used, up to now, in Middle East.
- "Ghubari" numerals, which are the ancestors of the modern numerals, were often used, during the Middle Ages, in Maghreb (North Africa) and Spain, at the time when the Arabic Muslim civilization was flourishing.

It is well established that the principle of local value was used by the Babylonians much earlier than by the Hindus [Cantor 1907] and that the Maya of Central America used the principle and symbols for zero in a well-developed numeral system of their own [Bowditch 1910], [Morley 1915]. We must recall that our study concerns the origin of our numerals (the "Ghubari" numerals) and not the Indian mathematics history. Professor Florian Cajori had well situated the problematic: ''The controversy on the origin of our numerals does not involve the question of the first use of local value and symbols for zero; it concerns itself only with the time and place of the first application of local value to the decimal scale and with the origin of the forms or shapes of our ten numerals'' [Cajori 1919].

Some assumptions have been put forward as to the origin of the numeral symbols. One of these assumptions, based essentially on the testimony of medieval Arabic writers from the Eastern Islamic world, and due to M. F. Woepcke [Woepcke 1863a] [Woepcke 1863b] considers the modern numerals as deriving from Indian characters. R. G. Kaye puts in doubt the testimonies of medieval Arabic writers from the Eastern Islamic world and try to invalidate the numeral Hindu origin hypothesis and to find an European origin of our numerals [Kaye 1907], [Kaye 1908], [Kaye 1911a], [Kaye 1911b], [Kaye 1918], [Kaye 1919]. In previous work [Boucenna 2006] we have given a new identification to the symbols of the "Ghubari" numerals through the mixed pagination of an Arabian Algerian manuscript, "Kitab khalil bni Ishak El Maliki", of the beginning of the 19th century. The correspondence, without ambiguity, between the modern numerals and their eldest the "Ghubari" numerals has been established. The elements of the "Abjadi" calculation ("Hissab El-joummel" or "Guematria"), particularly the notions of the numeric values ("Abjadi" numerical values) of the Arabic and Hebrew letters, have also been recalled. Then we showed the relation between the "Ghubari" numeral and the Arabic letter whose numerical value is equal to this numeral [table 1] and we gave the initial transformations imposed to the Arabic letters to become "Ghubari" numerals.

**Table 1**: The "Ghubari" numerals used in the numbering of the sheets
of the manuscript "Kitab khalil bni Ishak El Maliki" [Boucenna 2006]

| Numerical value | Arabic Letter | Likely Initial version of the Mashriki numeral | Ghubari numeral | European version of "Ghubar" Numeral | Final version of the Mashriki numeral | Hebrew letter |
|---|---|---|---|---|---|---|
| 1 | ا | ا | ا | ا | ١ | א |
| 2 | ب | ب | ح | ح | ٢ | ב |
| 3 | ج | ج | ج | ج | ٣ | ג |
| 4 | ح | ح | ح | 4 | ٤ | ד |
| 5 | ه | ه | 4 | ح | ٥ | ה |
| 6 | و | و | 6 | 6 | ٦ | ו |
| 7 | ز | ز | ) | ) | ٧ | ז |
| 8 | ح | ح | 8 | 8 | ٨ | ח |
| 9 | ط | ط | 9 | 9 | ٩ | ט |
| 10 | ي | ي |  |  | ٠ | י |
| 90 | ص |  | ص | O |  |  |

We also mentioned the permutation of the numbers 4 and 5 in Europe. In old European versions one recovers the order 4 and 5 used in Maghreb and in Spain [Hill 1999].

## 2. Problematic of the numeral origin

The history of the origin of the numerals resembles somewhat to the one of a child whose father is presumed to be one of the two men A and B. The testimonies written of several writers corroborate with the hypothesis that the child is the son of the A man. The AND tests show that he is rather the child of the B man. The B man personal objects have not reveal the existence of any relation between this man and the child's mother. Moreover, the B man nearest parents who testified that the child is the son of A man. The DNA test is given by this relation that exists between the shape of the ghubari numeral (and the Mashriki numeral) and the Arabic letter whose numerical value is equal to this numeral. This cannot be a fortuitous coincidence when it concerns, in an obvious way, 80% of the cases. The fact to say that our ghubari numerals are in fact Arabic letters is in agreement with the comparison, made by an unknown author, between the Western Arabic forms of the numerals and the letters of the Western Arabic script [Kunitzsch 2006]. Professor Paul Kunitzsch [Kunitzsch 2003], [Kunitzsch 2005], [Kunitzsch 2006] uses the term ''huruf al-ghubar'' taken from an Arabic reference that speaks well about "huruf" which means "letter". The Arabo-Islamic writers were well conscious that the numerals were letters, among the Arabic alphabet letters, as we have just shown [Boucenna 2006]. Professor Paul Kunitzsch suggested that new theories can only be brought forward when new, unknown, source material becomes available [Kunitzsch 2006] and emphasizes the need for further research to establish the existence of Hindu-Arabic numerals in the Western-Islamic sources prior to 1300 AD [Kunitzsch 2003].

The main argument of the partisans of the hypothesis of the Indian origin of the numbers is essentially based on the testimony of medieval Arabic writers from the Eastern Islamic world. Professor Florian Cajori [Cajori 1919] underlines the size of the difficulty to which was confronted R. G Kaye in his tentative to establish the hypothesis of a non-Hindu origin of the numerals. R. G. Kaye showed that this testimony is not without faille. A re-lecture of the testimony of medieval Arabic writers from the Eastern Islamic world impose itself. Al-Biruni testimony reported par Rabbi ben Ezra is at the origin of the majority of the other testimonies.

### 3. Kayes's lecture of testimony of Arabic medieval writers

#### 3.1. Review of medieval Arabic writers testimonies

R. G Kaye reviewed the testimonies of medieval Arabic writers from the Eastern Islamic world [Cajori 1919]. Severus Sebokht (662 AD) indicated that in the latter half of the seventh century the nine numerals were known in Arab lands and were attributed to the Hindus [Nau 1910], [Smith 1917]. Hurt by the alleged arrogance of certain Greek scholars, Sebokht praises the science of the Hindus and speaks of "their valuable methods of computation"... I wish only to say that this computation is done by means of nine signs. Unfortunately, he leaves it to us to guess whether or not he used the zero. About two centuries after Sebokht, appeared the famous arithmetic of the Arab Al- Khowarizmi (850 AD). The Arabic original manuscript is lost, but a Latin translation has come down to us under the title "Algoritmi de numéro Indorum". While this title refers to Indian numerals, they are not actually used in the book. A book on the astronomical tables of Al-Khowarizmi that was written by Muhammed ibn Ahmed Al-Biruni (973-1038) was translated into Hebrew by Rabbi ben Ezra, who says in his introduction that a Hindu astronomical work had been translated into Arabic and that, after the time of Alchowarizmi, "scholars do their multiplications, divisions, and extraction of roots as is written in the book of the [Hindu] scholar which they possess in translation" [Smith 1918]. Other Arabic authors who in the titles of their texts refer to the Hindus are enumerated by Kaye [Kaye 1911]. About 987 A.D. appeared "The great Treatise on the Table relating to the Indian Calculus". Soon after came "The Principles of the Indian Calculus". About 1030, "The satisfactory Treatise on Indian Arithmetic". There were two works, both bearing the same title, "Indian Arithmetic" one of the ninth century, the other of the tenth. A Latin text, attributed to Abraham, a Jew of whom little is known, is entitled "Liber augmenti et diminutionis vocatus numeratio divinationis, ex eo quod sapientes Indi posuerunt". The Latin authors took the place of the Arabic authors. The Italian Leonardo of Pisa, after travelling in Egypt, Syria, Greece, Sicily, wrote in 1202 his Liber abbaci in which he calls our numerals with the zero "figuras Indorum". The Byzantine morik, Maximus Planudes (1260-1330), wrote "Arithmetic according to the Hindus".

The raised number of the testimonies didn't convince R. G. Kaye. The evidence from these and some other texts that we have omitted, in favour of the Hindu origin of our numerals, is not so strong as one might think. In some cases no Hindu symbols are actually employed by the authors; the arithmetic and algebra set forth do not seem to bear Hindu characteristics. Kaye suspects that the word "Indian" was often incorrectly applied. Yet this testimony, as a whole, comes with a force that is difficult to break [Cajori 1919]. R. G. Kaye invokes other arguments to show the non Indian origin of the numerals [Cajori 1919]. If R G Kaye had information concerning the relation between the symbol of a numeral and the Arabic letter whose numerical value is equal to this numeral, he would have very quickly concludes.

**2.1. Linguistic confusion**

R G Kaye [Kaye 1907], [Kaye 1908], [Kaye 1911a], [Kaye 1911b], [Kaye 1918], [Kaye 1919] considered that a historic mistake occurred because of a confusion between the words "hindus", "Hindi" and "hindasi" in the translation of the writings of the Arabic authors. I. Taylor and M. F. Woepcke, and their followers, have ascribed to the Hindus the use of mathematical processes in early centuries, when, as a matter of fact, there is no evidence whatever to show that the Hindus actually used these processes at so early a date. This historical error arose according to Kaye in the mistranslation of the word "hindasi" [Cajori 1919]. R. G. Kaye was right to raise the possible role of a bad use of the words "hindus", "Hindi", and "hindasi" recovered in the testimony of medieval Arabic writers from the Eastern Islamic world and in the translations.

Let us reconsider the senses of the words "Hindi", and "hindasi". In Arabic language, the word "hind" has several significances. "Hind" is a woman's first name; "Hind" is also the name of an Arabic tribe of Arabia. "El Hind" is India. The word "hind" designates also "hundred camels" and/or "hundred camels and more" [Ibnou Mandhour a]. The word "hind" is bound thus to the notion of a number, big enough, without any relation with India.

In the Arabic sociology qualifying "Hindi" does not mean the Indian origin necessarily, it can mean: magic, bizarre, surprising and exotic. The fruits called in France "figues de barbari" (of the old Maghreb) are called in their own country, Maghreb, "Hindi", certainly because of their exotic appearance. "Hindi" means also: makes well, work of art, as "Nasl Hindi" (Hindi sword), or "Nasl Yamani" (Yemeni sword) that means: sword made well (without being achieved necessarily in India or in Yemen) [Ibnou Mandhour a]. To designate the "Geometry" the arabo-Islamic used the word "Hindassa". In fact, initially, "Hindassa" designated the techniques of the hydrological engineers (Mouhandiss) requiring the skills that are a matter for the genius. Their predictive power of the existence of the water underground expanses is astonishing. By extrapolation, all work of prediction and realization of imposing artwork is assigned to a "Mouhandiss". Thus, the word "Hindassi" takes the significance of: genius, brilliant. [Ibnou Mandhour b]. From this clarification, a direct relation appears between the words "Hindi" and "hindassi" in the sens that the two words can mean brilliant and genius. A brilliant or ingenious work is said " Hindassi " or " Hindi ".

Woepcke admits that ordinarily the word "hindasi" significates "geometrical", "measure", but asserts that this interpretation seemed impossible when used in connection with an explanation of the rule of "double false position" and the process of "casting out the nines" for the reason that these processes are purely arithmetical in nature [Woepcke 1863]. Because of the resemblance between the words "hindasi" and "hindi" or "Indian" Woepcke concluded that with the particular authors in question "hindasi" meant "Indian" and that, therefore, the "double false position" and "casting out the nines" were known to the early Hindus. This would seem to imply the use of our notation [Cajori 1919]. M. F. Woepcke did take into account only the "geometric" sense of the word "hindassi" and excluded the larger sense "ingenious" of the words "Hindi" and "hindassi" which can be applied to the purely arithmetical processes. Thus "Al-Hissab-Al-Hindi" can be translated simply by "Ingenious Calculation" and not by " Indian Calculation ".

**3. Al-Biruni Testimony**

The Testimony of Al-Biruni reported by Rabbi ben Ezra is at the origin of the majority of the other testimonies. Al-Biruni (973 -1050), a brilliant scholar and mathematician of the 11st century, has done several journeys in India and achieved a remarkable survey on the Indian civilization under all its aspects: customs, languages, sciences and geography. When Al-Biruni speaks about India he speaks as a knowledgeable. One of the important sources of information which we have about Indian origin of the numerals comes from Al-Biruni. Referring to the Indian numerals in a famous book written about 1030, Al-Biruni wrote: "**Whilst we use letters for calculation according to their numerical value, the Indians do not use letters at all for arithmetic. And just as the shape of the letters that they use for writing is different in different regions of their country, so the numerical symbols vary**". This testimony of Al-Biruni is the more considered by the supporter of the Indian origin of the numbers. Al-Biruni, well known mathematician of the arabo-Islamic civilization used the local value to the decimal scale and the forms or shapes of our ten numerals. Al-Biruni says well ''**Whilst we use letters for calculation according to their numerical value**''. One is right to ask: which would be these letters, having numerical values that were used by Al-Biruni, as numerals, to represent his numbers and to do his calculations? There is only one answer: the numerals used by Al-Biruni are well the Arabic letters having numerical values. This is in agreement with what we showed in previous work, where we affirmed that: The symbol of a numeral corresponds to the Arabic letter whose numerical value is equal to this numeral [Boucenna 2006]. For Al-Biruni, it does not make any doubt that he knows well that the numerals he uses are a lot of Arabic letters. Therefore, the testimony of Al-Biruni was not well understood.

In the second part of his testimony, Al-Biruni says: ''**the Indians do not use letters at all for arithmetic. And just as the shape of the letters that they use for writing is different in different regions of their country, so the numerical symbols vary**'' shows indeed that the numerals used by the Indian in all India, and which vary from an Indian region to another, have no relation with the numerals that Al-Biruni used.

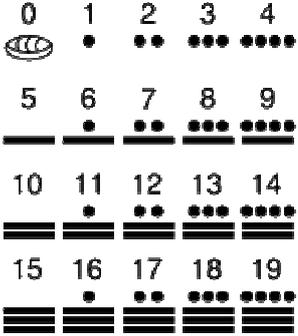

Fig. 1.  Mayan Numerals

Fig. 2.  Indian Numerals

This means that numerals used by the Indian in all India have no relation with our numerals (used by Al-Biruni). If Al-Biruni, had seen the numerals of the Maya (Fig. 1), he would have had the same astonishment that the one that he had while seeing the Indian numerals (Fig. 2). Contrary to the Greeks, the Romans, the Arabo-Islamic, the Hebrew ….who used their letters to represent the numbers (in the arithmetic), the Indians did not use letters at all for arithmetic. This proves that the Indian never used the numerals that Al-Biruni used and that we also use now. This testimony of the 11st century, in agreement with the conclusion of R G Kaye [Cajori 1919], shows that "the child's mother" did not have a relation with presumed "child's father".

## 4. Probable Origin of this confusion

Rabbi ben Ezra which translated A book on the astronomical tables of Al-Khowarizmi that was written by Muhammed ibn Ahmed Al-Biruni, says in his introduction that a Hindu astronomical work had been translated into Arabic and that, after the time of Alchowarizmi, "scholars do their multiplications, divisions, and extraction of roots as is written in the book of the [Hindu] scholar which they possess in translation [Smith 1918]. Al-Biruni would have spoken of an Indian origin of the calculation methods that would have been confused with the origin of our numerals. A bad interpretation of the testimony of Al-Biruni, in contradiction with his observation and his precedent testimony would be at the origin of the important number of all others testimonies of the medieval Arabic writers from the Eastern Islamic world and thereafter of the testimonies of the Latin authors and testimonies of the Arabic writers from the Western Islamic world.

All these testimonies are at the origin of the point of view of M. F. Woepcke and the modern western historians. The modern arabo-Islamic historians' opinion on the origin of the numerals follows that of the European historians [Kunitzsch 2005]. Carra de Vaux. [Vaux 1917] says that Al-Biruni must have drawn his information about Indian numerals from legend on Creation written in 943 A.D. due to the Arabic historian Massoudi (943) [Vaux 1911]. The report of the historian Massoudi concerning the origin of the numerals, could probably be founded on a rumour that would have been voluntarily propagated, at the time of the writting of the "Mashriki" numerals, for political considerations related to the history of that time. Indeed, it is surprising to see the existence of two systems in the Arabo-Islamic world : The "Ghubari" numerals in the west and the "Mashrikhi" numerals in the East. The " Ghubari " numerals would have taken their definitive shapes in Maghreb or in Spain, and would have passed to the Mashrek to be put back in shape to give the numerals "Mashriki".

Let us replace the invention of the numeral in its historic context. The Abassides (750-1258) have just decapitated the Omeyyades empire of Damascus (661 - 750) that consisted of the Maghreb (North Africa), Spain and the Mashrek (Middle East) until the valley of the Indus and the mountains of the Pamir [Petit 1913]. In 750, the Omeyyades are slaughtered and the Abbasside dynasty seized power. The Abbassides founded a new capital, Baghdad (762). The Maghreb stopped obeying to the Abbassides. The Omeyyade Abderrahmane, escaped the massacre and became master of the Omeyyade kingdom in Spain (755). The revolts of the kharidjite in Maghreb gave birth to the Roustomides kingdom in Tiaret, Algeria (776-909). Idris, proclaimed himself Caliph and Imam in Morocco, founded the dynasty of the Idrissides (788-985). The Aghlabides, with the consent of the Caliph Abbasside, Haroun Er-Rachid, governed Tunisia (800-909). El Ma'mûn (813-833) founded Beit El Hikma in Baghdad. El Ma'mûn (813-833) recruited a lot of skilled men and among them the mathematician Al-Khawarizmi. At that time (750-800), the Maghreb and Spain separate from the Abasside Empire of Baghdad. A big discovery took place: The Ghoubari numerals are born in Maghreb or in Spain. The discovery is too important to be ignored. But how to accept and to recognize that some very interesting things can take place in these rebel countries: Maghreb and Spain. The politics mingled, El Ma'mun ordered to his employers of Beit El Hikma to redraw other numerals from the Ghubari numerals and to invent them an origin to cancel their origin (rebel). India, the exotic and miraculous place, is the indicated country. Then a beautiful history was invented, the history of the trailer who was even received by the Caliph himself, and of the book, the book of the big Brahma. However, the ghubari numeral, with its direct relation with the Arabic letters and their numerical values continue to be used in Maghreb and in Spain and it passed to Europe to become our modern numeral. A meticulous study of the "Mashriki" numeral shows that his creator is not really the inventor of the numeral, he is an artist, of great knowledge, but obeys to the master's orders.

This does not exclude the hypothesis that probably the idea of the numerals could have had an embryonic state during the Omeyyade time in the Middle East or before. The ghubari numerals were really represented by letters in their initial version. The idea of the numeral would have emigrated to Maghreb or to Spain, probably with the Kharidjite, to finally give the "Ghubari" numerals ancestors of our modern numerals by the transformation of the Arabic letters. I believe that it is probably in the Roustomide Kingdom in Tiaret (Tihert, Algeria) that the "Ghubari" numerals would have been put in shape. My argument is based on the use, until now, of the Arabic letters by the Mozabites, heirs of the Roustomide Kingdom, to represent the numbers. The rules of calculation are written while using the equivalence Arabic letters - numbers. It is not in contradiction with the report of Severus Sebokht (662 AD) who speaks of methods of calculation only using new signs [Nau 1910, Smith 1917]. These signs are simply the first nine Arabic letters of the Abjadi order. He does not say if the numbers were represented using our numeral "zero" or, instead, the parenthesis or something else are used as with the Babylonians. Besides, even though the word "Hindi" is used in the sense of "Indian" and not in the sense of "ingenious" or "brilliant", Severus Sebokht speaks of "methods of calculation" achieved while using nine symbols, he does not speak of the representation of the numbers by our nine numerals.

## 8. CONCLUSION

The origin of our numerals that we inherited from the arabo-Islamic civilization remained one enigma. The hypothesis of the Indian origin remains, with controversies, without serious rival. It was the dominant hypothesis for more of one century. Its partisans found and constructed a lot of arguments. The testimonies of the medieval authors have been interpreted to its advantage. The opposite opinions have been dismissed and ignored. A voluntary amalgam between the history of our modern numerals and the Indian mathematics history is made. The rational contradictions often passed under silence.

A meticulous observation of the numerals allows us to affirm that our numerals are in fact more or less modified Arabic letters. The numerals do not have a simple resemblance with some Arabic letters, but every number looks like the Arabic letter whose numerical value is equal to this numeral. The element of the "Abjadi" calculation gives us a theoretical support, independent of the letters and numerals, which explains our observation. Besides a re-lecture of the testimonies of the medieval authors, particularly the testimony of Al-Biruni, that is probably at the origin of all others testimonies speaking of the Indian origin of the numerals, is in perfected agreement with the fact that our numerals are Arabic letters.

We have then a second way concerning the origin of our modern numerals which is only at its beginnings. Thorough researches are necessary to better understand the history of our numerals. A rigorous new reading of the medieval testimonies with a new mind imposes itself.